\documentclass[12pt]{article}
\usepackage{amsmath,amssymb,amscd,latexsym,euscript,mathrsfs,amsthm}
\usepackage[matrix,arrow,curve]{xy}

\newcommand{\Lim}{\underleftarrow{\lim}}
\newcommand{\Ab}{{\rm Ab}}

\DeclareMathOperator{\Ob}{Ob}
\DeclareMathOperator{\Mor}{Mor}

\DeclareMathOperator{\Ker}{Ker}

\DeclareMathOperator{\cd}{cd}
\DeclareMathOperator{\gld}{gl\,dim}
\DeclareMathOperator{\Dim}{Dim}
\DeclareMathOperator{\Imm}{Im}
\DeclareMathOperator{\fF}{\mathfrak{F}}

\newcommand{\NN}{{\,\mathbb N}}
\newcommand{\mC}{{\mathscr C}}
\newcommand{\mD}{{\mathscr D}}
\newcommand{\ZZ}{{\,\mathbb Z}}
\newcommand{\mA}{{\mathcal A}}

\newtheorem{theorem}{\bf Theorem}[section]
\newtheorem{lemma}[theorem]{\bf Lemma}
\newtheorem{proposition}[theorem]{\bf Proposition}
\newtheorem{corollary}[theorem]{\bf Corollary}
\newtheorem{definition}{\sc Definition}[section]
\newtheorem{example}[definition]{\sc Example}

\def\leq{\leqslant}
\def\geq{\geqslant}

\title
{
GLOBAL DIMENSION OF POLYNOMIAL RINGS 
IN PARTIALLY COMMUTING VARIABLES
}
\author{Ahmet A. Husainov}

\begin{document}

\maketitle

\begin{abstract}
For any free partially commutative monoid 
$M(E,I)$, we compute 
the global dimension of the category of $M(E,I)$-objects 
in an Abelian category with exact coproducts. 
As a corollary, we 
generalize Hilbert's Syzygy Theorem to polynomial rings  
in partially commuting variables. 
\end{abstract}



\section*{Introduction}

In this paper, the global dimension of the category of objects in an Abelian 
category with the action of free partially commutative monoid is 
computed. As a corollary, a formula for the global dimension 
of polynomial rings in partially commuting variables 
is obtained.

Let $\mA$ be any Abelian category. 
By \cite[Chapter XII, $\S 4$]{mac1963}, 
extension groups $Ext^n(A,B)$ are consisted 
of congruence classes of exact sequences 
$0\rightarrow B\rightarrow C_1\rightarrow \cdots \rightarrow C_n\rightarrow A$ 
in $\mA$ for $n\geq 1$  and $Ext^0(A,B)=Hom(A,B)$.
It allows us to define 
the {\em global dimension of $\mA$} by 
$$
	\gld \mA= \sup\{ n\in \NN: (\exists A,B\in \Ob\mA)~ Ext^n(A,B)\not=0\}.
$$
Here $\NN$ is the set of nonnegative integers.
(We set $\sup\emptyset =-1$ and $\sup\NN=\infty$.)
For a ring $R$ with $1$, 
$\gld R$ 
is the global dimension of the category of 
left $R$-modules.

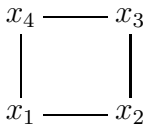
\begin{figure}[h]
$$
  \xymatrix
 {
	x_4	\ar@{-}[d]\ar@{-}[r] & x_3 \ar@{-}[d]\\
	x_1 \ar@{-}[r]	&  x_2
  }
$$
\caption{Pairs of commuting variables}\label{cycl4}
\end{figure}

As it is well known
 \cite[Theorem 4.3.7]{wei1994}, for any ring  
 $R$ with $1$, 
$$
	\gld R[x_1, \ldots, x_n]= n + \gld R.
$$
Moreover,  
by \cite[Theorem 2.1]{mit1978}, 
if  
$\mA$ is any Abelian category with exact coproducts and $\mC$ a {\em bridge category},   
then
$\gld \mA^{\mC}= 1 + \gld\mA$. It follows that   
$\gld \mA^{\NN^n}= n+ \gld \mA$ for the free commutative monoid $\NN^n$
generated by $n$ elements. 
We will get one of possible generalizations of this formula.
Let $M(E,I)$ be a free partially commutative monoid 
with a set of variables $E$, where  $I\subseteq E\times E$ 
is an irreflexive symmetric relation assigning the pairs 
of commuting variables. In this paper, we prove  that 
$$
\gld \mA^{M(E,I)}= n+ \gld \mA
$$
for any Abelian category with exact coproducts 
where  
$n$ is the sup of numbers
of mutually commuting distinct elements of $E$. 
For example, if $R[M(E,I)]$ is the polynomial ring   
in variables 
$E=\{x_1, x_2, x_3, x_4\}$ with the commuting pairs
$(x_i,x_j)$ corresponding to adjacent vertices of the graph demonstrated 
in Figure \ref{cycl4}, then for any ring  $R$ with $1$
 we have
$\gld R[M(E,I)] = 2+ \gld R$.

The free partially commutative monoids have numerous applications 
in combinatorics and computer sciences \cite{die1997}.
Our interest in their homology groups is concerned with the studying 
a topology of mathematical models for concurrency 
\cite{X20082}.

\section{Cohomology of small categories}

Throughout this paper let 
$\Ab$ the category 
of Abelian  groups and homomorphisms,
$\ZZ$ the additive group of integers, and $\NN$ the set of nonnegative
integers or the free monoid with only one generator. 
For any category  $\mA$ and a pair
 $A_1,A_2\in \Ob\mA$, denote by $\mA(A_1,A_2)$ the set of all morphisms
$A_1\rightarrow A_2$. 
A {\em diagram} $\mC\rightarrow \mA$ is a functor from a small category $\mC$ 
to a category $\mA$.
Given a small category $\mC$ we denote by
$\mA^\mC$  the category of diagrams  $\mC\rightarrow \mA$ and natural 
transformations.
For $A\in \Ob\mA$, let $\Delta_\mC A$ (shortly $\Delta A$) denote a diagram
$\mC \rightarrow \mA$ with constant values $A$ on objects and  
$1_A$ on morphisms. 
 
In this section, we recall some results from the cohomology theory 
of small categories.

\subsection{Homology groups of a nerve}

Recall a definition of a nerve of the category and properties of 
homology groups of simplicial sets. We refer the reader to  \cite{mac1963} and 
\cite{gab1967} for the proofs. 

\subsubsection
{\bf A nerve of the category} Let $\mC$ be a small category. Its {\em nerve} 
$N_*\mC$
is the simplicial set in which $N_n \mC$ 
consists of all sequences of composable morphisms 
$c_0 \stackrel{\alpha_1}\rightarrow c_1 \stackrel{\alpha_2}\rightarrow
\cdots \stackrel{\alpha_n}\rightarrow c_n$ in $\mC$ for $n>0$ and
$N_0\mC= \Ob\mC$. 
For $n>0$ and $0\leq i\leq n$, 
boundary operators $d^n_i: N_n\mC \rightarrow N_{n-1}\mC$
acts as 
$$
d^n_i(c_0 \stackrel{\alpha_1}\rightarrow 
\cdots \stackrel{\alpha_n}\rightarrow c_n)=
c_0 \stackrel{\alpha_1}\rightarrow 
\cdots 
\stackrel{\alpha_i}\rightarrow \hat{c_i} \stackrel{\alpha_{i+1}}\rightarrow
\cdots  \stackrel{\alpha_n}\rightarrow c_n ~.
$$
Here $c_0 \stackrel{\alpha_1}\rightarrow c_1 
\stackrel{\alpha_2}\rightarrow \cdots 
\stackrel{\alpha_i}\rightarrow \hat{c_i} \stackrel{\alpha_{i+1}}\rightarrow
\cdots  \stackrel{\alpha_n}\rightarrow c_n \in N_{n-1}\mC$ 
is the $(n-1)$-fold sequence 
obtained from 
$c_0 \stackrel{\alpha_1}\rightarrow c_1 \stackrel{\alpha_2}\rightarrow
\cdots \stackrel{\alpha_n}\rightarrow c_n$ for $0<i<n$ by
substitution the morphisms 
$c_{i-1} \stackrel{\alpha_i}\rightarrow c_i 
\stackrel{\alpha_{i+1}} \rightarrow c_{i+1}$ 
by their composition
$c_{i-1} \stackrel{\alpha_{i+1}\circ\alpha_i}\longrightarrow c_{i+1}$.
The map $d^n_0$ removes  $\alpha_1$ with $c_0$ and    
$d^n_n$ removes  $\alpha_n$ with $c_n$.
Degeneracy operators
$s^n_i: N_n\mC \rightarrow N_{n+1}\mC$ insert in
$c_0 \rightarrow \cdots \rightarrow c_n$ the identity morphism 
 $c_i \rightarrow c_i$ for every
$0 \leq i \leq n$. 

\subsubsection
{\bf Homology groups of simplicial sets}
Let $X$ be a simplicial set given by  boundary operators $d^n_i$ 
and degeneracy operators
$s^n_i$ for $0\leq i\leq n$. Consider a chain complex $C_*(X)$ of free Abelian groups
 $C_n(X)$ generated by the sets $X_n$ for $n\geq 0$.
Differentials
$d_n: C_n(X)\rightarrow C_{n-1}(X)$ are defined  on the basis elements $x\in X_n$
 by  $d_n(x)= \sum_{i=0}^n(-1)^i d^n_i(x)$. 
Let $C_n(X)=0$ for $n<0$. The groups $H_n(X)= \Ker d_n/ \Imm d_{n+1}$ 
are called  {\em $n$-th homology groups of the simplicial set
 $X$}.
The groups $H_n(X)$ are isomorphic to $n$-th singular homology groups of 
the geometric realization of $X$ by the Eilenberg theorem
 \cite[Appl. 2]{gab1967}.

\subsubsection
{\bf Cohomology of a category with coefficients in an Abelian group}
For a small category $\mC$, let $H_n(\mC)$ denote the $n$-th homology group 
of the nerve $N_*\mC$.
 For a simplicial set $X$ and an Abelian group  $A$, 
{\em cohomology groups} $H^n(X,A)$ are defined as cohomology groups of the 
complex  $Hom(C_*(X),A)$. Let $\mC$ be a small category. 
We introduce its  {\em cohomology groups} $H^n(\mC,A)$ 
with coefficients in $A$ as $H^n(N_*\mC,A)$. 
It follows from \cite[Chapter III, Theorem 4.1]{mac1963} that there 
is the following exact sequnce (Universal Coefficient Theorem)
$$
0 \rightarrow Ext(H_{n-1}(\mC), A)\rightarrow H^n(\mC,A)
\rightarrow Hom(H_n(\mC),A) \rightarrow 0
$$

\subsection{Cohomology of categories with coefficients in diagrams}

Recall the definition and properties of right derived functors 
$\Lim^n_{\mC}: \Ab^{\mC}\rightarrow \Ab$ of the limit functor. 

\subsubsection
{\bf Definition of cohomology of categories with coefficients 
in diagrams}
Let $\mC$ be a small category.
For every family 
$\{A_i\}_{i\in I}$ of Abelian groups we consider the 
direct product
$\prod\limits_{i\in I}A_i$ 
 as the Abelian group of maps 
$\varphi: I \longrightarrow \bigcup_{i \in I} A_i~$
such those $\varphi(i)\in A_i$ for all 
$i\in I$.

For any functor 
$F: \mC \rightarrow \Ab$, consider the sequnce of 
Abelian groups
$$
C^0({\mC},F) = 
\prod\limits_{c_0 \in \Ob\mC} F(c_0) , ~\ldots ,~
 C^n({\mC},F) = 
\prod\limits_{c_0 \stackrel{\alpha_1}\rightarrow \cdots 
\stackrel{\alpha_n}\rightarrow c_n}
F(c_n),~ \ldots
$$
and homomorphisms
$\delta^n: C^n (\mC,F) \rightarrow C^{n+1} (\mC,F)$ 
defined by
\begin{multline*}
  (\delta^n \varphi) (c_0 \stackrel{\alpha_1}\rightarrow \cdots
  \stackrel{\alpha_{n+1}}\rightarrow c_{n+1}) =\\
\sum_{i=0}^n (-1)^i \varphi (c_0 \stackrel{\alpha_1}\rightarrow \cdots
\stackrel{\alpha_i}\rightarrow {\hat c}_i
\stackrel{\alpha_{i+1}}\rightarrow \cdots
\stackrel{\alpha_{n+1}}\rightarrow c_{n+1}) +\\
(-1)^{n+1} F(c_n \stackrel{\alpha_{n+1}}\rightarrow
c_{n+1} ) ( \varphi(c_0 \stackrel{\alpha_1}\rightarrow \cdots
\stackrel{\alpha_n}\rightarrow c_n ) ).
\end{multline*}
Let $C^n(\mC,F)=0$ for $n<0$. 
The equalities $\delta^{n+1}\delta^n=0$ hold for all integer $n$. 
The obtained cochain complex will be denoted by
$C^*(\mC,F)$. Abelian groups
$H^n(C^*(\mC, F))= \Ker\delta^{n+1}/\Imm\delta^n$ are called  
{\em cohomology groups of the small category $\mC$ with  
coefficients in a diagram} $F$ and denoted by $\Lim^n_{\mC}F$.

It follows from
\cite[Appl. 2, Prop. 3.3]{gab1967} by the substitution  
 $\mA=\Ab^{op}$ that  
the functors $\Lim^n_{\mC}$ are   
$n$-th right satellites of 
$\Lim_\mC: \Ab^{\mC}\rightarrow \Ab$. 
Since the category $\Ab^{\mC}$ 
has enough injectives, the functors $\Lim^n_{\mC}$ 
are isomorphic to right derived of the limit functor.

\subsubsection
{\bf Cohomology of categories without retractions}
A morphism $\alpha: a \rightarrow b$ категории $\mC$ is a 
{\em retraction} if there exists a morphism 
 $\beta: b\rightarrow a$ such that $\alpha\beta= 1_b$.

\begin{proposition}{\rm \cite[Prop. 2.2]{X1997}}\label{normalization}
If a small category $\mC$ does not contain nonidentity retractions, 
then for any diagram $F:\mC\rightarrow \Ab$, the groups $\Lim^n_{\mC}F$
 are isomorphic to the homology groups of the 
subcomplex   $C^*_+({\mC},F)\subseteq C^*(\mC,F)$ composed of the products
\begin{displaymath}
 C^n_+({\mC},F) = \prod_{c_0 \stackrel{\not=}\rightarrow 
\cdots \stackrel{\not=}\rightarrow c_n}
F(c_n), \quad n \geq 0,
\end{displaymath}
where indices run the sequences  
$c_0 \stackrel{\alpha_1}\longrightarrow c_1 \stackrel{\alpha_2}\longrightarrow
\cdots \stackrel{\alpha_{n}}\longrightarrow c_{n}$ such those
$\alpha_i\not=id_{c_i}$ for all $1\leq i\leq n$.
\end{proposition}
\begin{corollary}\label{dlength}
If a small category $\mC$ does not contain nonidentity retractions 
and the length $m$ of every sequence of nonidentity morphisms
$c_0\stackrel{\alpha_1}\longrightarrow\cdots
\stackrel{\alpha_m}\longrightarrow  
c_m$ is not greater than $n$, then  $\Lim^k_{\mC}=0$ for $k>n$.
\end{corollary}

\begin{example}\label{parhom}
Let $\Theta$ be the category with $\Ob\Theta=\{a, b\}$ and 
$\Mor\Theta=\{1_a, 1_b,  a\stackrel{\alpha_1}\rightarrow b, 
a\stackrel{\alpha_2}\rightarrow b \}$. It follows from 
Proposition \ref{normalization}
that for any diagram   
$F: \Theta\rightarrow \Ab$ and $n>1$, the groups $\Lim^n_{\Theta}F$ equal $0$.
\end{example}

For any Abelian group $A$, $\Lim^n_{\mC}\Delta A\cong H^n(\mC,A)$.
\begin{lemma}\label{cdtors}
Let  
$\Theta^n$
 be the $n$-th power of the category $\Theta$ for $n\geq 1$. 
The functors $\Lim^k_{\Theta^n}$ equal $0$ for all $k>n$.
For any Abelian group $A$, there is an isomorphism 
$\Lim^n_{\Theta^n}\Delta A\cong A$.
\end{lemma}
\begin{proof}
The first assertion follows from Corollary \ref{dlength}.
Since the geometric realization of the nerve of $\Theta^n$ 
is the $n$-dimensional torus,  
$H_k(\Theta^n)\cong \ZZ^{n\choose k}$ for $0\leq k\leq n$. Here ${n\choose k}$ 
is the binomial coefficients. Universal Coefficient Theorem 
for the cohomology groups of the nerve of 
 $\Theta^n$ gives
 $H^n(\Theta^n,A)\cong A$. 
\end{proof}

\subsubsection
{\bf Strongly coinitial functors}
A small category $\mC$ is {\em acyclic} if  
$H_n(\mC)=0$ for all $n>0$ and $H_0(\mC)=\ZZ$.
Let  $S: \mC \rightarrow \mD$ be a functor from 
a small category to an arbitrary category.
For any $d\in \Ob\mD$, a  {\em fibre} (or {\em comma-category})
$S/d$ is the category which objects are given by pairs $(c, \alpha)$ 
where 
$c\in Ob(\mC)$ and $\alpha\in \mD(S(c),d)$.
Morphisms $(c_1, \alpha_1)\rightarrow (c_2, \alpha_2)$ in $S/d$ 
are triples $(f, \alpha_1, \alpha_2)$ with
$f\in \mC(c_1,c_2)$
satisfying $\alpha_2\circ S(f)=\alpha_1$.
If $S$ is a full embedding $\mC\subseteq \mD$, then $S/d$
is denoted by $\mC/d$.

\begin{definition}
A functor $S: \mC \rightarrow \mD$ between small categories
is called  {\em strongly coinitial} if $S/d$ is acyclic for each 
 $d\in \mD$.
\end{definition} 
\begin{lemma}[Oberst]\label{oberst}
 Let $\mC$ and $\mD$ be small categories. 
If $S:\mC \rightarrow \mD$ be a strongly 
coinitial functor, then the canonical homomorphisms 
$\Lim^n_{\mD}F\rightarrow \Lim^n_{\mC}FS$
are isomorphisms for all $n\geq 0$.
\end{lemma}
\begin{proof}
It follows from the opposite 
assertion  \cite[Теорема 2.3]{obe19682} for the functors 
$S^{op}: \mC^{op}\rightarrow \mD^{op}$ and
$F^{op}: \mD^{op}\rightarrow \Ab^{op}$.
\end{proof}

\subsection{Cohomological dimension of a small category}

Let  $\NN$ be the set of nonnegative integer numbers.
We will be consider it as the subset of 
$\{-1\}\cup\NN\cup\{\infty\}$ ordered by 
$-1 < 0 < 1 < 2 < \cdots < \infty$.

\begin{definition}
Cohomological dimension $\cd\mC$ of a small category $\mC$ is 
the sup in $\{-1\}\cup\NN\cup\{\infty\}$ 
of the set $n\in \NN$ for which the functors 
$\Lim^n_\mC: \Ab^\mC\rightarrow \Ab$ are not equal $0$.
\end{definition}
It follows from Lemma \ref{cdtors} that $\cd\Theta^n=n$. 
Lemma \ref{oberst} gives the following 
\begin{corollary}\label{hdcoinit}
If there exists a strongly coinitial functor
$S:\mC\rightarrow \mD$ between small categories, then
  $\cd\mC\geq \cd\mD$.
\end{corollary}

A subcategory $\mD\subseteq \mC$ is said to be {\em closed} 
if $\mD$ is a full subcategory
containing the domain for any morphism whose codomain is in $\mD$.

\begin{corollary}\label{fcoinit}
Let $\mD_j\subseteq \mD$ be a family of closed subcategories for 
all  $j\in J$.
If the inclusion 
$\bigcup_{j\in J}\mD_j\subseteq \mD$ is strongly coinitial, then
$\cd\mD = \sup_{j\in J}\{\cd\mD_j\}$.
\end{corollary}
\begin{proof}
 For $c\in \Ob\mC$, let $\mC_c\subseteq\mC$ 
be denote a full subcategory which consists of  
$c'\in \Ob\mC$ having morphisms
 $c'\rightarrow c$. It follows from \cite[Corollary 7]{mit1981} that
 the equality $\cd\mC= \sup_{c\in \Ob\mC}\cd\mC_c$ holds.
Consequently 
$\sup_{j\in J}\{\cd\mD_j\}= \cd \bigcup_{j\in J}\mD_j \leq \cd\mD$. 
Since  the inclusion  
 $\bigcup_{j\in J}\mD_j\subseteq \mD$ is strongly coinitial, the equality 
follows from Corollary \ref{hdcoinit}.
\end{proof}

\section{Dimension of a free partially 
commutative monoids}

We will prove the main results. 
We compute the Baues-Wirsching dimension of a 
free partially 
commutative monoids and
show a formula for the global dimension 
of the category of objects with actions 
of a free partially 
commutative monoid. We prove that 
for any graded  $R[M(E,I)]$-module, there exists a free resolution.

\subsection{Cohomological dimension of 
the factorization category}

We consider the category of factorization of a small 
category, although we applicate it 
for the case of the small category is a monoid.

\subsubsection
{\bf The category of factorizations}
Let $\mC$ be a small category. Objects of the {\em category of factorizations} 
$\fF\mC$ \cite{bau1985} are all morphisms of $\mC$. For any  
$\alpha, \beta \in \Ob\fF\mC= {\Mor}\mC$, the set of 
morphisms $\alpha\rightarrow \beta$ consists from all pairs 
$(f,g)$ of morphisms in $\mC$ satisfying  
 $g\circ\alpha\circ f = \beta$. The composition
of $\alpha\stackrel{(f_1,g_1)}\longrightarrow \beta$ and 
$\beta\stackrel{(f_2,g_2)}\longrightarrow \gamma$
 is defined by
$\alpha\stackrel{(f_1\circ f_2,g_2\circ g_1)}\longrightarrow \gamma$.
The identity of an object $a\stackrel\alpha\rightarrow b$ of 
$\fF\mC$ equals 
$\alpha \stackrel{(1_a, 1_b)}\longrightarrow \alpha$.

\subsubsection
{\bf Baues-Wirsching dimension}
A {\em natural system of Abelian groups on $\mC$} is any 
functor  $F: \fF\mC\rightarrow \Ab$. 
Baues and Wirsching introduce 
cohomology groups  $H^n(\mC,F)$ of $\mC$ with coefficients in 
a natural system $F$ and have proved that these groups are isomorphic 
to $\Lim^n_{\fF\mC}F$. 
The {\em Baues-Wirsching dimension} $\Dim \mC$ is the cohomological 
dimension of  $\fF\mC$.

\begin{example}
Let $\NN=\{1, a, a^2, \ldots \}$ be the free monoid generated 
by one element. It easy to see that the inclusion 
$\Theta_a \subseteq \fF\NN$ of the full subcategory with the objects 
$\Ob\Theta_a = \{1,a\}$ is strongly coinitial. The subcategory  
$T_a$ is closed in $\fF\NN$. It is isomorphic to $\Theta$ from 
Example  \ref{parhom}. Consequently $\Dim\NN= 1$. 
\end{example}

\begin{proposition}\label{degr}
For any integer $n\geq 1$,  
$\Dim\NN^n= n$. 
\end{proposition}
\begin{proof}
Consider the full subcategory $\Theta_a^n \subseteq 
\fF\NN^n$ with objects $(a^{\varepsilon_1}, \ldots, a^{\varepsilon_n})$
where $\varepsilon_i\in \{0,1\}$ for all $1\leq i\leq n$. 
It not hard to see that it is isomorphic to $\Theta^n$ and the 
fibre of the inclusion over 
$(a^{k_1}, \ldots, a^{k_n})\in \NN^n$ is isomorphic to the product
$\Theta_a/a^{k_1}\times \cdots \Theta_a/a^{k_n}$. 
Since $H_i(\Theta_a/a^k)=0$ for  
$i>0$ and  $H_0(\Theta_a/a^k)\cong \ZZ$, it follows that 
the category 
$\fF\NN^n$ contains the srongly coinitial subcategory 
$\Theta_a^n$, which is isomorphic to 
$\Theta^n$. It is clear that  $\Theta_a^n$ is closed in $\fF\NN^n$.
Hence, $\Dim\NN^n= \cd\Theta^n= n$. 
\end{proof}

\subsection{The dimension of a free partially commutative monoid}

This subsection is devoted to computing the Baues-Wirsching dimension 
of free partially commutative monoids.

\subsubsection
{\bf The independence graph}
Let $E$ be a set and $I\subseteq E\times E$ an 
irreflexive symmetric binary relation on $E$.
Monoid given by a generating set  $E$ and relations  
$ab=ba$ for all $(a,b)\in I$ is called  
{\em free partially commutative} and denoted by 
$M(E,I)$.
\begin{figure}[h]
$$
\xymatrix{
	 b \ar@{-}[d]\ar@{-}[r] & a \ar@{-}[d] \ar@{-}[rd]  \\
	 c \ar@{-}[r] & d \ar@{-}[r] & e
}
$$
\caption{The independence graph}\label{grph}
\end{figure}
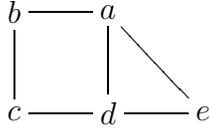

The pair $(E,I)$ may be considered as a simple 
{\em independence graph of $M(E,I)$}
with the set of vertices $E$ and edges $\{a, b\}$ for all pairs 
 $(a,b)\in I$. 
It is shown in Figure \ref{grph} the independence graph 
of the monoid given by the generators 
$E=\{a,b,c,d,e\}$ and relations 
$ab=ba$, $bc=cb$, $cd=dc$, $ad=da$, $ae=ea$, $de=ed$.

The {\em clique number $\omega(E,I)$} of a simple graph 
with verices  
$E$ and edges $I$ is the sup of cardinalities 
of its finite complete subgraphs. If $(E,I)$ contains complete graphs $K_n$ 
for all $n\in \NN$, then $\omega(E,I)=\infty$.

For example, the clique number 
of the graph in Figure \ref{grph} is equal to $3$.

\subsubsection
{\bf Computing the dimensions of 
free partially commutative monoids}
Let $V$ be the set of maximal cliques of the independence graph
of $M(E,I)$. (These cliques may be infinite.) 
For example, the set $V$ for the graph in Figure \ref{grph}, 
consists of the sets 
 $\{a, b\}$, 
$\{b, c\}$, $\{c, d\}$, $\{a, d, e\}$.
Let $E_v$ be the set of verices belonging to a clique $v$ and  
$M(E_v)$ the submonoid of $M(E,I)$ generated by 
$E_v$. It is clear that $M(E_v)$ are commutative monoids.
The category of factorization $\fF M(E_v)$ is a closed subcategory of 
 $\fF M(E,I)$.

\begin{lemma}\label{coinit}{\rm\cite{X20083}}
The inclusion $\bigcup\limits_{v\in V}\fF M(E_v)\subseteq \fF M(E,I)$ 
is strongly coinitial.
\end{lemma}

\begin{theorem}
$\Dim M(E,I)= \omega(E,I)$.
\end{theorem}
\begin{proof}
For every subset of mutually commuting elements
$\{e_1, \ldots, e_n\}\subseteq E$, the full subcategory
 $\fF M(\{e_1, \ldots, e_n\})$ 
is closed in $\fF M(E,I)$. Hence, the equality is true in the case of
$\omega(E,I)= \infty$.
The subcategories $\fF M(E_v)$ are closed in $\fF M(E,I)$. 
It follows from Lemma \ref{coinit} 
that we can use Corollary \ref{fcoinit}.
We get $\cd \fF M(E,I)= \sup_{v\in V}\{\cd\fF M(E_v)\}$.
If $E_v$ are finite, then we get the assertion $\Dim M(E,I)= \omega(E,I)$  
by Proposition  \ref{degr}.
\end{proof}

\subsection{The generalized syzygy theorem} 

In this subsection, we prove the main theorem.

\subsubsection
{\bf The global dimension of the category of $M(E,I)$-objects}
By 
\cite[Corollary 13.4']{mit1972} 
for any small category $\mC$ 
and Abelian category with exact coproducts, 
there exists the inequality 
$\gld \mA^{\mC}\leq \dim\mC+ \gld \mA$. Here $\dim$ is the 
{\em Hochschild-Mitchell dimension}. We will show that 
if $\mC=M(E,I)$, then the equality holds. 
\begin{theorem}\label{main1}
Let $\mA$ be an Abelian category with exact coproducts.
Then 
$\gld \mA^{M(E,I)}= \omega(E,I)+ \gld \mA$.
\end{theorem}
\begin{proof}
For any finite subset $E'\subseteq E$, the submonoid generated by 
 $E'$ is cancellative \cite{die1997}. It follows that 
 $M(E,I)$ is cancellative and  $\Dim M(E,I)$ is equal to Hochschild-Mitchell
 dimension 
$\dim M(E,I)$ \cite[Theorem 3.1]{X1997}.
Consequently 
$\gld \mA^{M(E,I)}\leq \omega(E,I)+ \gld \mA$.
For each finite subset 
of mutually commuting elements
 $S\subseteq E$ 
there exists a retraction $M(E,I)\rightarrow M(S)$. 
It follows by  
\cite[Corollary 1.4]{mit1978} that 
$\gld \mA^{M(E,I)}\geq \gld \mA^{M(S)}= |S| + \gld \mA$. 
Since $\dim M(E,I)$ is equal to sup of cardinalities 
 $|S|$ of finite subsets $S\subseteq E$ of mutually commuting elements, 
we get 
$\gld \mA^{M(E,I)}\geq \omega(E,I)+ \gld\mA$.
\end{proof}

\subsubsection
{\bf Graded syzygies}
Let  $R$ be a ring with с $1$.  
The monoid ring has the natural graduation 
$R[M(E,I)]= \bigoplus\limits_{n\in \ZZ}R[M(E,I)]_n$ 
by  $R$-modules 
$R[M(E,I)]_n=\{r\mu: r\in R, \mu\in M(E,I), 
|\mu|=n\}$. In particular $R[M(E,I)]_0=R$. Let $R[M(E,I)]_n=0$ for all $n<0$.
The ring $R$ with $1$ is called {\em projective free} if 
any projective $R$-module is free.
By \cite[$\S$8.7, Corollary 2]{bur1987} and Theorem \ref{main1} we 
get:
\begin{corollary}
Let $M(E,I)$ be a free partially commutative monoid and
 $R$ projective free ring with $\gld R= n < \infty$. 
If there is the maximal number $m < \infty$ of mutually commuting  distinct  
elements of $E$, then for each  bounded below $\ZZ$-graded
$R[M(E,I)]$-module $A$, there exists an exact sequence 
of $\ZZ$-graded $R[M(E,I)]$-modules and $\ZZ$-graded homomorphisms 
of degree $0$
$$
0\rightarrow F_{n+m} \rightarrow F_{n+m-1} \rightarrow \cdots 
\rightarrow F_0 \rightarrow A \rightarrow 0\,,
$$
with free bounded below $\ZZ$-graded $R[M(E,I)]$-modules 
$F_0$, $F_1$, ... , $F_{n+m}$.
\end{corollary}

\end{document}